\documentclass[11pt]{amsart}

\usepackage{amsxtra}

\newtheorem{theorem}{Theorem}[section]
\newtheorem{proposition}[theorem]{Proposition}
\newtheorem{lemma}[theorem]{Lemma}
\newtheorem{corollary}[theorem]{Corollary}

\theoremstyle{definition}

\theoremstyle{remark}

\newtheorem{remark}[theorem]{Remark}

\numberwithin{equation}{section}

\newcommand{\D}{\mathcal{D}}
\renewcommand{\epsilon}{\varepsilon}

\newcommand{\loc}{{\rm loc}}

\newcommand{\R}{\mathbb{R}}

\DeclareMathOperator{\supp}{supp}
\DeclareMathOperator{\tr}{tr}

\title[Lieb-Thirring inequalities]{Lieb-Thirring inequalities on
  the half-line with critical exponent}
 
\author{Tomas Ekholm}
\address{Tomas Ekholm, Complexo Interdisciplinar, Av. Prof. Gama Pinto
  2, 1649-003 Lisboa, Portugal}
\email{tomase@math.kth.se}

\author{Rupert L. Frank}
\address{Rupert L. Frank, Royal Institute of
  Technology, Department of Mathematics, 100 44 Stockholm, Sweden}
\email{rupert@math.kth.se}

\begin{document}

\begin{abstract}
We consider the operator $- \frac {d^2}{dr^2} - V$ in
$L_2(\R_+)$ with Dirichlet boundary condition at the origin. For the moments
of its negative eigenvalues we prove the bound
\begin{equation*} 
\tr \left(- \frac {d^2}{dr^2} - V\right)_-^\gamma \leq C_{\gamma, \alpha}
\int_{\R_+} \left( V(r) - \frac {1}{4r^2} \right)_+^{\gamma + \frac {1
  + \alpha}2} r^\alpha \, dr
\end{equation*}
for any $\alpha \in [0,1)$ and $\gamma + \frac {1 + \alpha}2 \geq 1$.
\end{abstract}

\maketitle


\section{Introduction}

In this paper we consider inequalities for moments of negative
eigenvalues of one-dimensional Schr\"odinger operators. If $V$ is a
real-valued function on $\R$ which vanishes at infinity (at least in
some averaged sense) then the negative spectrum of
$-\frac{d^2}{dx^2}-V$ consists of discrete eigenvalues of finite
multiplicities. The celebrated Lieb-Thirring inequality states that 
\begin{equation}\label{eq:lt}
\tr_{L_2(\R)}\left(-\frac{d^2}{dx^2}-V\right)_-^\gamma \leq L_{\gamma}
\int_{\R} V(x)_+^{\gamma+\frac12} \,dx 
\end{equation}
holds with a constant $L_{\gamma}$ independent of $V$ if and only if
$\gamma \geq 1/2$. Here and in the sequel $v_\pm:=\max\{\pm v,0\}$
denotes the positive or negative part of $v$. In the non-critical
case $\gamma>1/2$ bound \eqref{eq:lt} was proved in \cite{LiTh},
and the sharp values of the constants $L_\gamma$ for $\gamma\geq 3/2$
were found in \cite{LiTh} and \cite{AiLi}. The inequality in the
endpoint case $\gamma = 1/2$ was established much later by Weidl
\cite{Wei}. In this case, Hundertmark-Lieb-Thomas \cite{HuLiTh}
determined the sharp value of the constant $L_{1/2}$. The sharp
constants for $1/2 < \gamma < 3/2$ are still unknown. 

Egorov-Kondrat'ev \cite{EgKo} studied weighted versions of
inequality \eqref{eq:lt}. For any $\alpha > 0$ they show that
\begin{align}\label{eq:ek}
  \tr_{L_2(\R)}\left(-\frac{d^2}{dx^2} - V\right)_-^\gamma
  \leq C_{\gamma,\alpha}^{EK} 
  \int_\R V(x)_+^{\gamma + \frac{1+\alpha}2} |x|^\alpha \,dx
\end{align}
holds if and only if $\gamma > \frac{1 + \alpha}2$. Note that the 
endpoint case is excluded in contrast to the case $\alpha = 0$ and that
the potential $V$ can only appear with a power strictly larger than
one. Of course, the location of the origin in \eqref{eq:ek} is
arbitrary and can be optimized over. 

Let us turn to half-line Schr\"odinger operators $-\frac{d^2}{dr^2}-V$
in $L_2(\R_+)$ where $\R_+=(0,\infty)$. Throughout we impose a
Dirichlet boundary condition at the origin. By the variational
principle, inequalities \eqref{eq:lt} and \eqref{eq:ek} remain valid
for these operators with the same constants, i.e., 
\begin{align} \label{eq:ekhalfline}
  \tr_{L_2(\R_+)}\left(-\frac{d^2}{dr^2} - V\right)_-^\gamma
  \leq C_{\gamma,\alpha}^{EK} \int_{\R_+} V(r)_+^{\gamma + \frac{1 +
  \alpha}2} r^\alpha \,dr. 
\end{align}
For the case $\alpha > 0$ it was shown in \cite{EgKo} that
\eqref{eq:ekhalfline} holds even for smaller values of $\gamma$ than
\eqref{eq:ek}, namely for all
$\gamma>\max\{\frac{1-\alpha}2,0\}$. However, the validity in the
endpoint case $\gamma=\frac{1-\alpha}2$ if $0\leq\alpha<1$ was left
open. 

In this paper we shall prove a substantially stronger inequality in
this endpoint case. Our main result, Theorem \ref{maintheorem}, says
that 
\begin{align} \label{eq:mainintro}
\tr_{L_2(\R_+)} \left(- \frac {d^2}{dr^2} - \frac 1{4r^2} -
V\right)_-^\gamma \leq C_{\gamma, \alpha} \int_{\R_+} V(r)_-^{\gamma +
  \frac {1 + \alpha}2} r^\alpha dr
\end{align}
for any $0\leq \alpha< 1$ and $\gamma\geq \frac{1-\alpha}2$. Since the
eigenvalues of $-\frac{d^2}{dr^2} - V$ decrease if we subtract
$\frac1{4r^2}$, inequality \eqref{eq:mainintro} extends
\eqref{eq:ekhalfline} to the endpoint case. Conversely, it is not
difficult to see that both \eqref{eq:ekhalfline} and 
\eqref{eq:mainintro} fail if $\gamma<\frac{1-\alpha}2$ and
$0\leq\alpha<1$, see Remark \ref{remarkoncriticalgamma}. 

The main advantage of \eqref{eq:mainintro} over \eqref{eq:ekhalfline} 
however lies in the Hardy term $\frac1{4r^2}$. First of
all, recall the (sharp) Hardy inequality 
\begin{equation}\label{eq:hardy}
\frac 14\int_{\R_+} \frac{|u|^2}{r^2}\,dr 
\leq \int_{\R_+} |u'|^2 \,dr,
\quad u\in C_0^\infty(\R_+).
\end{equation}
This implies that the operator $-\frac{d^2}{dr^2}-\frac1{4r^2}$ is
non-negative and that the constant $\frac14$ is maximal with respect
to this property. Now it is clarifying to rewrite
\eqref{eq:mainintro} as 
\begin{equation}\label{eq:mainintro2}
\tr \left(- \frac{d^2}{d r^2} - V \right)_-^\gamma
\leq C_{\gamma, \alpha} \int_{\R_+} \left( V(r) - \frac
1{4r^2}\right)^{\gamma +\frac{1+\alpha}2}_+ r^\alpha \, dr. 
\end{equation}
This is an inequality of the same form as \eqref{eq:ekhalfline} but
with a different right hand side. Inequality \eqref{eq:mainintro2}
shows that only the part of the potential which is larger than the
Hardy weight is necessary to control negative eigenvalues. In
particular, it follows from \eqref{eq:hardy} that $- \frac{d^2}{d r^2}
- V$ has no negative eigenvalues if $V(r)\leq \frac1{4r^2}$ for all
$r$. This is reflected in \eqref{eq:mainintro2} but not in
\eqref{eq:ekhalfline}. 

In the case $\alpha=0$ the right hand side of \eqref{eq:ekhalfline}
coincides, up to a constant, with the semi-classical phase-space
integral 
\begin{equation*}
\frac 1{2\pi} \iint_{\R\times\R_+} \!\!\left(k^2-V(r)\right)_-^\gamma
  \,dk\,dr. 
\end{equation*}
However, this semi-classical approximation does not take the repulsive
Di\-rich\-let condition at the origin into account. This is achieved by
\eqref{eq:mainintro2}, which decreases the relevant phase-space
integral considerably. Inequality \eqref{eq:mainintro2} can indeed be 
viewed as an infinite phase-space renormalization. 

Note also that the operator $- \frac {d^2}{dr^2} - \frac 1{4r^2}$
appears as the radial part of the \emph{two}-dimensional Laplacian
after the natural change of measure. Hence \eqref{eq:mainintro}
estimates moments of eigenvalues of the operator $-\Delta-V(|x|)$ in
$L_2(\R^2)$ corresponding to angular momentum zero. 

Our interest in inequality \eqref{eq:mainintro} originates partially from
our previous work \cite{EkFr}, where we proved a similar inequality in
the case $\alpha \geq 1$ and $\gamma > 0$. This was the main tool to
extend the multi-dimensional version of \eqref{eq:lt} in the same way
as \eqref{eq:mainintro} extends \eqref{eq:ekhalfline}. Note, however,
that in these considerations the endpoint case $\gamma=0$ is naturally
excluded. We mention also the recent alternative proof \cite{FrLiSei}
of the main result of \cite{EkFr}. 

The proof of \eqref{eq:mainintro} in the endpoint case
$\gamma = \frac{1-\alpha}2$ encounters several difficulties. The proof
in \cite{HuLiTh} of \eqref{eq:lt} for $\gamma = 1/2$ relies heavily
on the translation invariance of the whole-line operator. The earlier
proof of \cite{Wei} does so too, but to a lesser extent, and its
generalization to our non-translation invariant setting requires additional
ideas both on a conceptual and on a technical level. One crucial
ingredient in our proof is the combiniation of Neumann bracketing with
the ground-state representation. Despite this (certainly non-optimal)
approach we obtain reasonable values for the constant
$C_{\gamma,\alpha}$ in \eqref{eq:mainintro}. In the important special
case $\alpha = 0$, $\gamma = 1/2$ we work out upper and lower
bounds which differ by less than a factor $2.25$. 

In the final section of this paper we show how our main result can be
applied to yield a Lieb-Thirring inequality for the operator
associated with the generalized Hardy inequality 
\begin{align}\label{eq:hardygeneral}
\frac{(\sigma -1)^2}4  \int_{\R_+} \frac{|u|^2}{r^{2 - \sigma}} \, dr
\leq \int_{\R_+} r^\sigma |u'(r)|^2 \, dr
\end{align}
for suitable $u$. We mention in closing that inequality
\eqref{eq:mainintro} was useful when proving Lieb-Thirring
inequalities on regular metric trees \cite{EkFrKo}. 

{\it Acknowledgements.}
This work was partially supported by FCT Portugal, post-doc grant SFRH/BPD/23820/2005, (T.E.) and by the Swedish Foundation for International Cooperation in Research and Higher Education (STINT) (R.F.). R.F. would like to thank E.H. Lieb and R. Seiringer for hospitality at Princeton University and for helpful discussions. Remarks by A. Hansson are gratefully acknowledged.


\section{Main result}

Let $V:\R_+\to\R$ with $V_-\in L_{1,\loc}(\R_+)$ and $V_+\in
L_{p}(\R_+,r^\alpha dr)$ for some $\alpha\geq 0$ and some $1\leq
p<\infty$. The Schr\"odinger operator $- \frac {d^2}{d r^2} -
\frac{1}{4 r^2} - V$ in $L_2(\R_+)$ is defined via the closure of the
quadratic form 
\begin{equation*}
\int_{\R_+} \left(|u'|^2 - \frac{|u|^2}{4r^2} - V|u|^2 \right)\,dr
\end{equation*}
on $C_0^\infty(\R_+)$. (The fact that this form is bounded from below
will follow as soon as we have proved \eqref{eq:result} for, say, all
bounded $V$ with compact support.) Our main result is 

\begin{theorem} \label{maintheorem}
  Let $\gamma>0$ and $\alpha \in [0,1)$ such that $\gamma +\frac{1+\alpha}2 \geq 1$, then
  \begin{align} \label{eq:result}
    \tr \left(- \frac {d^2}{d r^2} - \frac{1}{4 r^2} - V \right)_-^\gamma
    \leq C_{\gamma,\alpha} 
    \int_{\R_+} V(r)^{\gamma + \frac {1+\alpha}2}_+ r^\alpha\, dr
  \end{align}
  with a constant $C_{\gamma, \alpha}$ independent of $V$. In the
  special case $\alpha = 0$, $\gamma = 1/2$ the sharp constant in
  this inequality satisfies $0.533 \leq C_{1/2,0} \leq
  1.185$.
\end{theorem}

Our result can also be stated in terms of the operator $-\frac {d^2}{d
  r^2} - V$, defined similarly as above with a Dirichlet boundary
condition at the origin. Then Theorem \ref{maintheorem} implies (see
\cite{EkFr} for a careful argument) 

\begin{corollary}
Let $\gamma>0$ and $\alpha \in [0,1)$ such that $\gamma
  +\frac{1+\alpha}2 \geq 1$, then 
\begin{align*}
\tr \left(- \frac {d^2}{d r^2} - V \right)_-^\gamma
\leq C_{\gamma, \alpha} 
\int_{\R_+} \left( V(r) - \frac 1{4r^2}\right)^{\gamma + \frac {1 +
    \alpha}2}_+ r^\alpha \, dr
\end{align*}
with the constant $C_{\gamma, \alpha}$ from \eqref{eq:result}.
\end{corollary}

\begin{remark}
  The most important estimate in Theorem \ref{maintheorem} is that for
  the critical case $\gamma = \frac{1-\alpha}2$ when $V$ appears with the
  exponent one on the right hand side of \eqref{eq:result}. It shows
  that eigenvalue moments of any order $0 < \gamma \leq 1/2$ can
  be estimated \emph{linearly} in $V$. (For scaling reasons however,
  the integral of $V$ now has to include a weight.) 
  This is in sharp contrast to the whole-line case \eqref{eq:lt} and
  \eqref{eq:ek}, where only moments of order $\gamma = 1/2$ can be
  estimated linearly, and where moreover the inclusion of a weight does 
  not allow for smaller values of $\gamma$. 
\end{remark}

\begin{remark}
  The operator $H_0:=- \frac {d^2}{d r^2}-\frac1{4r^2}$ has a virtual
  level, in the sense that $H_0-V$ has a negative eigenvalue for
  any non-negative $V\not\equiv 0$. This shows immediately that it is
  \emph{impossible} to estimate the number of negative eigenvalues of
  $H_0-V$ in terms of a (weighted) $L_p$-norm of $V$. In
  particular, the critical case $\gamma = 0$ is excluded in
  \eqref{eq:result} for $\alpha\geq 1$. In order to
  estimate eigenvalue moments of arbitrarily small order $\gamma > 0$ in
  terms of a (weighted) $L_p$-norm of $V$ it is necessary that the
  lowest eigenvalue $\lambda(\beta)$ of $H_0 - \beta V$ disappears
  faster than any polynomial as $\beta \to 0$+. Indeed, $\lambda(\beta)$
  is exponentially small in our case, see \cite{EkFr} for details.
\end{remark}

\begin{remark} \label{remarkoncriticalgamma}
  The condition $\gamma \geq \frac{1-\alpha}2$ in Theorem
  \ref{maintheorem} is sharp. Indeed, if  $\gamma < \frac{1-\alpha}2$
  then $V$ appears with a sublinear power in the right hand side of
  \eqref{eq:result}. Hence if we choose a sequence of potentials $V_n
  := n \chi_{(R,R+n^{-1})}$ with $R>0$ arbitrary, then the right hand
  side of \eqref{eq:result} tends to zero as $n\to\infty$. On the other hand,
  the sequence $- \frac {d^2}{d r^2} -\frac1{4r^2} - V_n$ converges in
  norm resolvent sense to the operator $- \frac {d^2}{d r^2}
  -\frac1{4r^2} - \delta_R$, which has a negative eigenvalue. (This can be
  proved along the lines of \cite[Theorem~3.2.3]{AGHH}.) Hence 
  the limit of the left hand side of \eqref{eq:result} is positive. 
\end{remark}

\begin{remark}
The bounds on $C_{1/2,0}$ are based on numerical evaluation of Bessel
functions. Note that the upper bound differs from the lower bound by
less than a factor $2.25$. It is remarkable that $C_{1/2,0}$ is
strictly larger than $1/2$, which is the sharp constant in \eqref{eq:lt} for
$d=1$ and $\gamma = 1/2$. This means that the (repulsive) Dirichlet
boundary condition at the origin cannot completely compensate the
(attractive) potential $\frac 1{4r^2}$. In particular we prove that a
potential well $V$ situated near a finite $R$ may have a lower ground-state 
energy than the same well translated to $R=\infty$. 
\end{remark}


\section{Proof of the Lieb-Thirring inequality}\label{sec:proof}

This section contains the proof of our main result, Theorem
\ref{maintheorem}. It will be given in Subsection \ref{sec:proofsub}
after we have stated two basic ingredients in Subsection
\ref{sec:interval}.

\subsection{Operators on a finite interval}\label{sec:interval}

Throughout this section we fix a constant $b > 0$. We define the
operator $H_b$ in $L_2(b,b+1)$ via the quadratic form 
\begin{equation*}
h_b[u] := \int_b^{b+1} \left| \frac {d}{dr} \left( \frac
    {u(r)}{\sqrt r} \right) \right|^2 r \, dr, 
\quad u\in H^1(b,b+1).
\end{equation*}
Note that this can also be written as
\begin{equation}\label{eq:hbbdry}
h_b[u] = \int_b^{b+1} \left(|u'|^2 - \frac{|u|^2}{4r^2}\right) \, dr 
- \frac{|u(b+1)|^2}{2(b + 1)} + \frac{|u(b)|^2}{2 b}.
\end{equation}
It follows that $H_b$ acts as $- \frac{d^2}{dr^2} -
\frac 1{4r^2}$ on functions satisfying the natural boundary conditions
\begin{align}\label{eq:bc}
u'(b) - \frac {u(b)} {2b } = u'(b+1) - \frac {u(b+1)}{2(b+1)} = 0.
\end{align}
As an aside we remark that $H_b$ coincides with the two-dimensional
Neumann Laplacian in $\{x \in \R^2 : b < |x| < b + 1\}$ restricted to
radially symmetric functions. 

For any $k>0$ the resolvent $(H_b + k^2)^{-1}$ exists and is an
integral operator with kernel $G_b(\cdot,\cdot,k)$, i.e.,
\begin{align}\label{eq:reskernel}
\left((H_b + k^2)^{-1}u \right)(r) = \int_b^{b+1} G_b(r,s,k)u(s)\, ds,
\quad r\in [b,b+1].
\end{align}
We shall need

\begin{lemma} \label{lemmaC}
  For any $b\leq r\leq b+1$ the function $G_b(r, r, \, \cdot \,) :
  \R_+ \to \R$ is continuous, non-negative and non-increasing. 
\end{lemma}

Using the explicit expression of $G_b$ in terms of Bessel
functions we shall establish

\begin{proposition} \label{propositionD}
  Let $k>0$ and $0\leq\alpha<1$. Then there is a constant
  $C_\alpha(k)>0$ such that for all $b>0$ and $b\leq r\leq b+1$ one has
  \begin{align}\label{eq:uniformbound}
    G_b(r, r, k) \leq C_\alpha(k) \, r^\alpha. 
  \end{align}
  For $\alpha=0$ and $k=3.555$ one may choose $C_0(3.555)=1/3$.
\end{proposition}

The proof of this proposition will be given in Subsection
\ref{sec:green} below. Now we use the result to estimate the
lowest eigenvalue of the Schr\"odinger operator $H_b-V$ on the
interval $(b,b+1)$.

\begin{corollary} \label{intevcor}
  Assume that $V$ satisfies
  \begin{align} \label{eq:vsmall1}
    \int_b^{b+1} V (r) r^\alpha\, dr \leq C_\alpha(k)^{-1}
  \end{align}
  for some $k>0$, some $0\leq\alpha<1$ and $C_\alpha(k)$ from 
  \eqref{eq:uniformbound}. Then the lowest eigenvalue
  $\lambda$ of the operator $H_b - V$ satisfies
  \begin{equation*}\label{eq:intevcor}
  	-\lambda \leq k^2.
  \end{equation*}
\end{corollary}	

\begin{proof}
  By a standard approximation argument we may assume that $V$ is
  continuous. Moreover, by the variational principle we can restrict
  ourselves to the case $V \geq 0$ and $V\not\equiv 0$. Taking
  $u(r)=\sqrt r$ as a trial function we see that $\lambda<0$. We
  denote by $N(t^2)$ the number of eigenvalues of $H_b - V$ less than
  $-t^2$. By the Birman-Schwinger principle (see, e.g., \cite{BiSo2})
  we have
  \begin{align*}
    N (t^2)
    \leq \tr \sqrt V (H_b + t^2)^{-1} \sqrt V
    = \int_b^{b+1} V(r) G_b (r,r,t) \, dr.
  \end{align*}
  (Here we used the continuity of $V$ for the evaluation of the
  trace.) Now we let $t^2 \to -\lambda$ from above and use monotone
  convergence for the right hand side. Denoting
  \begin{equation*}
    F_b(t) := \max_{b\leq r \leq b+1} r^{-\alpha} G_b (r,r,t)
  \end{equation*}
  we find
  \begin{align*}
    1 \leq \int_b^{b+1} V(r) G_b (r,r,\sqrt{-\lambda}) \, dr
    \leq F_{b} (\sqrt{-\lambda}) \int_b^{b+1} V(r) r^\alpha\, dr.
  \end{align*}
  Combining this with the assumption \eqref{eq:vsmall1} we arrive at
  \begin{equation*}
    C_\alpha(k) \leq F_{b} (\sqrt{-\lambda}).
  \end{equation*}
  On the other hand, Proposition \ref{propositionD} is equivalent to
  $F_{b} (k) \leq C_\alpha(k)$. Since $F_{b}$ is a non-increasing
  function by Lemma \ref{lemmaC}, we conclude that
  $\sqrt{-\lambda}\leq k$, as claimed.
\end{proof}

The second ingredient in the proof of Theorem \ref{maintheorem} is the
following Poin\-ca\-r\'e-Sobolev inequality.

\begin{proposition} \label{propositionA}
  Let $0\leq\alpha<1$. Then there exists a constant $S_\alpha>0$ such
  that for all $b>0$ and for all $v \in H^1 (b,b+1)$ with
  $\int_b^{b+1} v(r) r \, dr =0$ one has
  \begin{align}\label{eq:propositionA}
    \max_{b \leq r \leq b+1} |v(r)|^2 r^{1-\alpha} 
    \leq S_\alpha \int_b^{b+1} |v'(r)|^2 r \, dr.
  \end{align}
  For $\alpha=0$ the sharp constant is $S_0 = 1/3$.
\end{proposition}

The proof of this proposition will be given in Subsection
\ref{sec:embedding} below. We remark that as $b$ grows
the function $\sqrt r$ on the interval $(b,b+1)$ becomes `almost
constant', so at least intuitively one recovers the inequality
\begin{align*}
  \max_{0 \leq x \leq 1} |v(x)|^2
  \leq \frac 13 \int_0^{1} |v'(x)|^2 \, dx,
  \qquad \int_0^{1} v(x) \, dx =0,
\end{align*}
which played an important role in \cite{Wei}. Note that allowing for
finite values of $b$ does not increase the constant.

Now we deduce from Proposition \ref{propositionA} an integral condition on $V$ that guarantees that the operator $H_b-V$ has only one negative eigenvalue.

\begin{corollary} \label{corollaryB}
  Assume that $V\not\equiv 0$ is a non-negative function on $(b,b+1)$
  satisfying
  \begin{align} \label{eq:vsmall}
    \int_b^{b+1} V (r) r^\alpha\, dr \leq S_\alpha^{-1}
  \end{align}
  for some $0\leq\alpha<1$ and $S_\alpha$ from \eqref{eq:propositionA}. Then the operator $H_b - V$ has exactly
  one negative eigenvalue.
\end{corollary}

\begin{proof}
  The existence of a negative eigenvalue has already been established
  in the proof of Corollary \ref{intevcor}. To prove the uniqueness we
  note that in view of Proposition \ref{propositionA} we have, for all
  $u \in H^1(b,b+1)$ with $\int_b^{b+1} u(r) \sqrt r \,dr = 0$, the
  inequality
  \begin{align*}
    h_b[u] - \int_b^{b+1} V (r) |u(r)|^2\, dr 
    \geq \left(1- S_\alpha \int_b^{b+1} V (r) r^\alpha \,dr\right)
    h_b[u].
  \end{align*}
  Since this is non-negative by \eqref{eq:vsmall}, we deduce by the
  variational principle that $H_b - V$ has at most one negative
  eigenvalue.
\end{proof}


\subsection{Proof of the main theorem}\label{sec:proofsub}

Throughout this section we fix $0\leq\alpha<1$. Our proof follows and
extends the ideas of \cite{Wei}. We divide it into four steps.

\emph{Step 1.}
It suffices to prove Theorem \ref{maintheorem} in the special case
$\gamma = \gamma_c:=\frac {1 - \alpha}2$. Indeed, the case $\gamma >
\gamma_c$ is already contained in \cite{EkFr} or, alternatively, may 
be deduced from the result for $\gamma = \gamma_c$ by the
argument of Aizenman-Lieb \cite{AiLi}. The latter is based on the
identity 
\begin{align*}
  B_{s,t} \lambda^s_- 
  = \int_{\R_+} \mu^{s-t-1}(\lambda+\mu)_-^t\,d\mu,
  \quad s > t,
\end{align*}
with some finite constant $B_{s,t}$ (which can be expressed in terms
of the beta function). Using it twice and assuming that the result is
proven in the critical case one obtains for any
$\gamma > \gamma_c$ that
\begin{align*}
& \tr\left(-\frac{d^2}{d r^2}-\frac{1}{4 r^2}-V\right)_-^\gamma \\
& \qquad \leq B_{\gamma,\gamma_c}^{-1} 
  \int_{\R_+} \mu^{\gamma - \gamma_c - 1} \tr \left(- \frac{d^2}{d
  r^2} - \frac 1{4 r^2} - V + \mu \right)_-^{\gamma_c} \, d\mu \\
& \qquad \leq B_{\gamma, \gamma_c}^{-1} C_{\gamma_c, \alpha}
  \int_{\R_+} \int_{\R_+} \mu^{\gamma - \gamma_c - 1} 
  (V(r) - \mu)_+^{\gamma_c + \frac {1 + \alpha}2} \, d\mu \,r^\alpha\,dr\\
& \qquad = B_{\gamma, \gamma_c}^{-1}
  B_{\gamma + \frac {1 + \alpha}2, \gamma_c + \frac {1 + \alpha}2}
  C_{\gamma_c, \alpha} \int_{\R_+} V(r)_+^{\gamma + \frac {1+\alpha}2}
  \, r^\alpha \, dr. 
\end{align*}
This is inequality \eqref{eq:result}, and so it remains to prove
the result for $\gamma = \gamma_c$.

\emph{Step 2.}
Now we begin with the main argument. The basic strategy is to divide
$\R_+$ into intervals such that the restriction of $- \frac {d^2}{d
  r^2} - \frac{1}{4 r^2} - V$ to these intervals has at most one
negative eigenvalue. The choice of boundary conditions for the
restricted operators is essential to achieve this. 
We choose boundary conditions \eqref{eq:bc} which come naturally with
the groundstate representation formula, see \eqref{eq:subst} below.

We may assume that $V\not\equiv 0$ is non-negative and, by standard
approximation arguments, that it has compact support in $\R_+$. Fix
$k>0$ arbitrary and let $\Psi(k) := \max \{S_\alpha, C_\alpha(k)\}$ where
$S_\alpha$, $C_\alpha(k)$ are the constants from Propositions
\ref{propositionA} and \ref{propositionD}. We set $a_1 := \min \supp
V$ and define a sequence $a_1 < a_2 < \ldots$ recursively by
\begin{equation}\label{eq:a_k}
  \int_{a_j}^{a_{j+1}} V (r) r^\alpha \, dr 
  = \frac 1{\Psi(k) (a_{j+1} - a_j)^{1-\alpha}}.
\end{equation}
This recursion stops when $a_N \geq \max\supp V$. The sequence is always 
finite since  $a_{j+1} - a_j \geq (\Psi(k) \|V\|_{L_1(r^\alpha
dr)})^{-1/(1-\alpha)} >0$, and clearly it covers $\supp V$. We set $a_0:=0$ 
and $a_{N+1} := \infty$.

Similarly as in the previous section we define operators $L_j$ in
$L_2(a_j,a_{j+1})$ via the quadratic form 
\begin{align*}
  \int_{a_j}^{a_{j+1}} \left| \frac {d}{dr} 
  \left( \frac{u(r)}{\sqrt r} \right) \right|^2 r \, dr
\end{align*}
with domain $H^1(a_j,a_{j+1})$ if $1\leq j\leq N$. If $j=0$ we consider the
closure of this form defined on $C_0^\infty(0,a_1]$.   
Note that for $u\in C_0^\infty(\R_+)$ one has
\begin{align}\label{eq:subst}
  \int_{\R_+} \left(|u'|^2-\frac{|u|^2}{4r^2}\right)\ dr =
  \int_{\R_+} \left| \frac {d}{dr} \left( \frac {u(r)}{\sqrt r}
  \right) \right|^2 r \, dr.
\end{align}
The variational principle implies that imposing natural boundary
conditions does not increase the operator, i.e.,
\begin{equation*}
  -\frac {d^2}{dr^2}-\frac{1}{4r^2} - V 
  \geq \bigoplus_{j=0}^N L_j - V.
\end{equation*}
Since $L_0-V=L_0\geq 0$ and similarily for $j=N$, we find that
\begin{align}\label{eq:bracketing}
  \tr \left(- \frac {d^2}{d r^2} - \frac{1}{4 r^2} - V \right)_-^{
    \frac {1-\alpha}2}
  \leq \sum_{1 \leq j \leq N-1} \tr (L_j - V)_-^{ \frac {1-\alpha}2}.
\end{align}

It remains to estimate $\tr (L_j - V)_-^{ \frac {1-\alpha}2}$ for
fixed $1 \leq j < N$. We shall implement a unitary change of variables
in order to obtain an operator on an interval of unit length and then
apply the results from Subsection \ref{sec:interval}. We put
\begin{align*}
  l_j := a_{j+1}-a_j, \quad b_j:=a_j/l_j,
\end{align*}
and introduce the unitary dilation operator $\mathcal U_j:
L_2(a_j,a_{j+1})\to L_2(b_j,b_j+1)$,
\begin{align*}
  (\mathcal U_j u)(r) := \sqrt{l_j} u(l_j r).
\end{align*}
One obtains the unitary equivalence
\begin{equation*}\label{eq:scaling}
  \mathcal U_j^{-1} (H_{b_j}-V_j) \mathcal U_j = l_j^2 (L_j-V)
\end{equation*}
where $V_j(r):=l_j^2 V(l_j r)$. Note that this potential satisfies
\begin{equation*}
  \int_{b_j}^{b_j+1} V_j(r) r^\alpha\,dr = \Psi(k)^{-1}
\end{equation*}
by \eqref{eq:a_k}. The definition of
$\Psi(k)$ together with Corollaries \ref{corollaryB} and
\ref{intevcor} implies that $H_{b_j}-V_j$ has exactly one negative
eigenvalue, and that its modulus does not exceed $k^2$. Combining this
with the above unitary equivalence and using \eqref{eq:a_k} once more we obtain
\begin{align*}
  \tr (L_j - V)_-^{ \frac {1-\alpha}2}
  \leq l_j^{-1 + \alpha} k^{1 - \alpha}
  =  k^{1 - \alpha} \Psi(k) \int_{a_j}^{a_{j+1}} V(r) r^\alpha \, dr. 
\end{align*}
In view of \eqref{eq:bracketing} this concludes the proof of
inequality \eqref{eq:result}.

\emph{Step 3.}
We next prove the upper bound $C_{1/2, 0} \leq 1.185$. For this
we note that the above proof yields 
\begin{equation*}
C_{\frac{1-\alpha}2,\alpha} 
\leq \inf_{k > 0} k^{1-\alpha} \max \left\{S_\alpha, C_\alpha(k)\right\}.
\end{equation*}
If $\alpha=0$ we choose $k=3.555$ and use
Propositions~\ref{propositionD} and \ref{propositionA} to get the
claimed estimate. See Remark \ref{conjck} concerning this choice.

\emph{Step 4.}
Finally, we prove the lower bound $C_{1/2, 0} \geq
0.533$. We shall first establish
\begin{equation}\label{eq:lowerbound}
	C_{\frac{1-\alpha}2,\alpha}
  \geq \sup_{R>0} R^{1-\alpha} I_0(R) K_0(R),
  \qquad 0\leq\alpha<1.
\end{equation}
For $\beta,R>0$ one can define the operator $- \frac {d^2}{d r^2} -
\frac{1}{4 r^2} - \beta \delta_R$ in a standard way via a quadratic
form. It follows from general principles that this operator has at
most one negative eigenvalue. Moreover, one easily establishes that
for any given $R$ there exists a unique $\beta=\beta(R)$ such that the
operator has $-1$ as an eigenvalue. Solving the eigenvalue equation
explicitly we obtain 
\begin{align*}
	u(r) = \left\{ 
  \begin{array}{ll}
  	\sqrt r \, I_0 (r) K_0 (R), & 0 < r < R, \\
    \sqrt r \, K_0 (r) I_0 (R), & R < r, 
  \end{array} \right.
\end{align*}
and simplifying with the help of the Wronski identity
\begin{equation}\label{eq:wronski}
	I_1(r)K_0(r) + I_0(r)K_1(r) = 1/r,
\end{equation}
see \cite[9.6.15]{AbSt}, we find
\begin{align*}
	\beta(R) 
	= \frac {u'(R-)-u'(R+)}{u(R)}
  = \frac 1{R\, I_0 (R) K_0 (R)}.
\end{align*}
By an approximation argument as in Remark~\ref{remarkoncriticalgamma}
one easily obtains the lower bound \eqref{eq:lowerbound}. 

Now assume that $\alpha=0$. Using the asymptotic behavior of the
Bessel functions one finds that $\beta(R) \to 2$ \emph{from below} as
$R \to \infty$, and hence $C_{1/2, 0} > 1/2$. Moreover, using
once again again Bessel function properties one can prove that $\beta$
has a unique minimum. Numerically one finds that it occurs at
$R=1.075$ and satisfies $\beta(1.075)^{-1}=0.533$. 

This concludes the proof of Theorem \ref{maintheorem}.


\section{The operators on a finite interval}

\subsection{Green's function}\label{sec:green}

By Sturm-Liouville theory (see, e.g., \cite{Wm}) we find that the resolvent
kernel \eqref{eq:reskernel} is given by
\begin{align}\label{eq:explicitkernel}
  G_b(r,s,k) = \left\{ 
  \begin{array}{ll}
    \frac {g_b(r, k) g_{b+1}(s, k)} {W_b(k)}, 
    & b \leq r \leq s \leq b+1, \cr  
    \frac {g_{b+1}(r, k) g_b(s, k)} {W_b(k)}, 
    & b \leq s \leq r \leq b+1, 
  \end{array}
  \right.
\end{align}
where
\begin{align*}
  g_c(r, k) & := 
  \sqrt r \left( I_1(c k) K_0 (k r) + K_1(c k) I_0(k r)\right),
  \quad c\in\{b,b+1\},\\
  W_b(k) & := I_1((b+1)k) K_1(bk) - I_1(bk) K_1((b+1)k).
\end{align*}
Here $I_n$ and $K_n$ denote the modified Bessel functions of the first
and second kinds of order $n$, see \cite{AbSt}.  

Now we give the simple

\begin{proof}[Proof of Lemma \ref{lemmaC}]
The continuity of $G_b(r,r,k)$ in $k$ follows from
\eqref{eq:explicitkernel} by the continuity of the Bessel
functions. Moreover, for any $f\in L_2(b,b+1)$ and any $k \geq t\geq
0$ one has 
\begin{align*}
0 \leq ((H_b + k^2)^{-1} f, f) \leq ((H_b + t^2)^{-1}f,f).
\end{align*}
Choosing $f$ as an approximate delta-function one easily finds that
the resolvent kernel is non-negative on the diagonal and
non-increasing in the spectral parameter. 
\end{proof}

Now we turn to the

\begin{proof}[Proof of Proposition \ref{propositionD}]
For fixed $k>0$ and $0\leq\alpha <1$ we define
\begin{align*}
  g_\alpha (x,b) & : = (b + x)^{-\alpha} G_b (b + x, b + x, k) \\
  & = (b+x)^{1 - \alpha} \big(I_1(bk) K_0((b+x)k) + I_0((b+x)k)
  K_1(bk)\big) \\ 
  & \qquad\times \big(I_1((b+1)k) K_0((b+x)k) + I_0((b+x)k) K_1((b+1)k)\big) \\
  & \qquad\times \big(I_1((b+1)k) K_1(bk) - I_1(bk) K_1((b+1)k)\big)^{-1}
\end{align*}
for $x \in [0,1]$, $b > 0$. We have to prove that there exists a
constant $C_\alpha(k) > 0$ such that for all $x \in [0,1]$, $b > 0$ one
has
\begin{align}\label{eq:uniformbound1}
  g_\alpha (x,b) \leq C_\alpha(k).
\end{align}

We begin with the case $\alpha=0$. 
Using the asymptotic behavior of the Bessel functions one finds
that, uniformly in $x\in[0,1]$,
\begin{align}\label{eq:bsmall}
  g_0(x,0) := \lim_{b \to 0} g_0(x,b) 
  &= \frac{x I_0(kx)(I_1(k)K_0(kx) + K_1(k)I_0(kx))}{I_1(k)}
\end{align}
and
\begin{align*}
  g_0(x,\infty) := \lim_{b \to \infty} g_0(x,b) 
  = \frac {\cosh (kx) \cosh (k(1-x))}{k \sinh k}.
\end{align*}
Both limiting functions $g_0(\cdot,0)$ and $g_0(\cdot,\infty)$ are
uniformly bounded on $[0,1]$. Since $g_0$ is continuous on
$[0,1]\times\R_+$ we obtain the bound \eqref{eq:uniformbound1} for
$\alpha=0$.

The statement is proved similarly for $0 < \alpha < 1$. Indeed, for
$b \geq 1$ the statement is weaker than for $\alpha = 0$. To treat small
$b$ one notices that $g_\alpha (x,b) \sim - (b + x)^{1 - \alpha} \log
(b + x)$ as $(x,b) \to (0,0)$.

Finally, we give a numerical estimate of $C_0(k)$ for the special choice
$k=3.555$. The function $g_0$ can be maximized numerically on
$[0,1]\times[0,\infty)$. (Simple estimates show that one can restrict
oneself to a compact subset.) One finds that the maximum is attained
at $(x,b)=(1,0)$ and one has 
\begin{equation*}
  \sup g_0 = g_0(1,0)= 0.333316 < 1/3.
\end{equation*}
Hence $C_0(3.555) < 1/3$.
\end{proof}

\begin{remark}\label{conjck}
  Numerical calculations suggest that the Green's function on the diagonal
  $G_b(r,r,k)$ attains its maximum at the right endpoint $r=b+1$ for any 
  value of $b$ and $k$. This would imply that
  \begin{align}\label{eq:conjck}
    \max_{x\in[0,1]} g_0(x,b) = g_0(1,b), \qquad b>0.
  \end{align}
  Moreover, one can check that the function $g_0(1,\cdot)$ attains its
  supremum in the limit $b\to 0$. By \eqref{eq:bsmall} and
  \eqref{eq:wronski} one finds the value $g(1,0)=\frac{I_0(k)} {k
  I_1(k)}$. Hence we believe that the sharp 
  constant in \eqref{eq:uniformbound} for $\alpha=0$ is given by
  \begin{align*} 
    C_0(k)=\frac{I_0(k)} {k I_1(k)}.
  \end{align*}
  Note that the RHS is a decreasing function and that $k=3.555$ in the
  above proof is chosen (almost) maximal with the property that
  $I_0(k)/k I_1(k)\leq1/3$. 
\end{remark}

\subsection{A Poincar\'e-Sobolev inequality}\label{sec:embedding}

We turn now to the proof of Proposition \ref{propositionA}. The core
is contained in

\begin{lemma}\label{minimization} 
Let $b > 0$ and $b \leq c \leq b + 1$. Then for all $v \in H^1(b, b +
1)$ satisfying $\int_b^{b+1} v(r) r \, dr 
= 0$ one has
\begin{align*}
|v(c)|^2 \leq \Phi(b,c) \int_b^{b+1} |v'(r)|^2 r \, dr
\end{align*}
with the sharp constant
\begin{align*}
\Phi (b,c) := \frac {1}{4(2b+1)^2} & \big( 4 (b+1)^4 \log(b+1) - 4 b^4
\log b \\
& - (2b+1)(3 + 6b + 6b^2 - 4c^2 + 4(2b^2 + 2b + 1) \log c )\big).
\end{align*}
\end{lemma}

\begin{proof} 
We shall assume $b < c < b+ 1$. The remaining cases $b=c$ and $b=c+1$
are proved similarly. We consider the functional 
\begin{align*}
F[v] := \frac {\int_b^{b+1} |v'(r)|^2 r \, dr} {|v(c)|^2}
\end{align*}
on the domain
\begin{align*}
\D [F] := \left\{ v \in H^1(b,b+1) :\, \int_b^{b+1} v(r) r \, dr = 0,\,
  v(c)\neq 0 \right\}.
\end{align*}
From the compactness of the embedding $H^1(b,b+1)\subset C[b,b+1]$, see, e.g., \cite{Ad}, it follows
that $F$ has a minimizer $v$. 
We may normalize $v$ by
\begin{align*}
\frac 1 {v(c)} \int_b^{b+1} v'(r)^2 r \, dr =1.
\end{align*}
In a standard way we derive the Euler-Lagrange equation $(v'(r)r)' =
\frac{2r}{2b+1}$ for $r \in (b,c) \cup (c,b+1)$ and the boundary
conditions $v'(r)=0$ for $r \in \{b,b+1\}$. We conclude that  
\begin{align*}
v(r) = \left\{ \begin{array}{ll}
D - \frac {b^2}{2b + 1} \log r + \frac {r^2}{2(2b+1)}, & \ b \leq r \leq c, \cr
D + \log c - \frac {(b+1)^2}{2b+1} \log r + \frac {r^2}{2(2b+1)}, & \
c \leq r \leq b+1,
\end{array}
\right.
\end{align*}
where
\begin{align*}
D := \frac {1}{4(2b+1)^2} \big( & 4 (b+1)^4 \log(b+1) - 4 b^4 \log b \\
& - (2b+1) (3 + 6b + 6b^2 - 2c^2 + 4(b+1)^2 \log c )\big).
\end{align*}
For the minimal value of the functional we obtain $F[v] =
v(c)^{-1}=\Phi(b,c)^{-1}$, as claimed. 
\end{proof}

\begin{proof}[Proof of Proposition \ref{propositionA}]
For $x\in[0,1]$, $b>0$, $0 \leq \alpha < 1$ we put 
\begin{align*}
\phi_\alpha(x,b) := & \, (b+x)^{1-\alpha}\Phi(b,b+x) \\
= & \, \frac {(b+x)^{1-\alpha}}{4(2b+1)^2} \, ( 4 (b+1)^4 \log(b+1) - 4 b^4
\log b \\
& \qquad \qquad \qquad - (1 + 2b)(3 + 6b + 6b^2 - 4(b + x)^2)) \\
& - \frac{2b^2 + 2b + 1}{2b+1}(b+x)^{1-\alpha}\log(b+x).
\end{align*}
We have to prove that there exists a constant $S_\alpha>0$ such that
for all $x\in[0,1]$, $b>0$ one has 
\begin{equation}\label{eq:sobolevconst}
	\phi_\alpha(x,b) \leq S_\alpha.
\end{equation}

First we note that $\phi_\alpha$ can be extended continuously to the boundary $\{b=0\}$. Indeed, uniformly in $x\in [0,1]$,
\begin{equation*}
\phi_\alpha(x,0) := \lim_{b\to 0} \phi_\alpha(x,b) = 
\left\{
\begin{array}{ll}
x^{1-\alpha}\left(-\log x +x^2-\frac34\right) & \quad \text{if}\ x\neq 0, \\
0 & \quad \text{if}\ x= 0.
\end{array}
\right.
\end{equation*}
Similarly, one finds that uniformly in $x\in [0,1]$,
\begin{equation*}
\phi_\alpha(x,\infty) := \lim_{b\to \infty} \phi_\alpha(x,b) = 
\left\{
\begin{array}{ll}
\frac 13 - x + x^2 & \quad \text{if}\ \alpha= 0, \\
0 & \quad \text{if}\ 0<\alpha <1.
\end{array}
\right.
\end{equation*}
Hence the continuity of $\phi_\alpha$ implies that
\eqref{eq:sobolevconst} holds with some finite constant $S_\alpha$. 

Finally, we turn to the issue of the sharp value of the constant for
$\alpha=0$. Numerically, one finds that $S_0=1/3$, which is
attained for $x\in\{0,1\}$ as $b\to\infty$. For the reader who feels
uncomfortable with this numerical optimization we sketch an
analytical proof below. We shall write $\phi$ instead of $\phi_0$. Note
that the limit of $\phi(\cdot,b)$ as $b\to\infty$ implies that the
sharp constant cannot be less than $1/3$. 

To prove the opposite inequality one checks first that $\phi(0,\cdot)$
is an increasing function with $\phi(0,0)=0$ and
$\phi(0,\infty)=1/3$. Similarly, $\phi(1,\cdot)$ is an
increasing function with $\phi(1,0)=1/4$ and
$\phi(1,\infty)=1/3$. Now we distinguish according to whether
$\phi(\cdot,b)$ has a local maximum in $(0,1)$ or not. In the latter
case we use the facts mentioned above to get 
\begin{align*}
\max_{x\in[0,1]}\phi(x,b) = \max\{\phi(0,b),\phi(1,b)\}\leq 1/3.
\end{align*}
Now consider the case where $x_0\in (0,1)$ is a local maximum. Again by the facts mentioned above it suffices to prove that
\begin{align}\label{eq:localmax}
	\phi(x_0,b)\leq 1/3.
\end{align}
We first claim that one necessarily has
\begin{equation}\label{eq:necessary}
	0 \leq x_0 \leq 1/\sqrt{6}
	\quad\text{and}\quad 
	0 \leq b \leq (1+\sqrt{5})/4.
\end{equation}
Indeed, note that $\partial_x^2 \phi_0 = \frac{1}{(2b+1)(b+x)}(6(b+x)^2-(2b^2+2b+1))$. Since $x_0$ is a local maximum, we conclude that
\begin{align}\label{eq:xestimate}
	0\leq x_0 \leq \sqrt{(1+2b+2b^2)/6}-b
	\quad\text{and}\quad
	\sqrt{(1+2b+2b^2)/6}\geq b,
\end{align}
which is easily seen to imply \eqref{eq:necessary}.

To proceed, we decompose $\phi(x,b)=\phi^{(1)}(x,b)+\phi^{(2)}(x,b)$, where
\begin{align*}
\phi^{(1)}(x,b) := \, & \frac {2b^2 + 2b + 1}{6(2b+1)}
\left(-6(b+x)\log (b+x) + 6(x + b - 1) \right. \\ 
& \qquad\qquad\qquad \left. + 3(x + b - 1)^2 - (x + b - 1)^3\right).
\end{align*}
We can estimate $\phi^{(1)}(x,b)\leq 0$ for all $(x,b)$ satisfying \eqref{eq:necessary} (with $x_0$ replaced by $x$).
Now we note that $\phi^{(2)}(\cdot,b)$ is a polynomial of degree three. A tedious but elementary calculation shows that it has a local maximum $x_1(b)$ and a local minimum $x_2(b)$ satisfying $x_1(b)<0<1/\sqrt6<1/2<x_2(b)$ for $0\leq b\leq (1+\sqrt 5)/4$. Hence we conclude that $\phi^{(2)}(x,b)\leq\phi^{(2)}(0,b)\leq 1/3$ for all $(x,b)$ satisfying \eqref{eq:necessary} (with $x_0$ replaced by $x$). This proves \eqref{eq:localmax}.
\end{proof}

\begin{remark} 
  What we actually have shown in the preceding proof is that the
  minimizer $v_*$ of the problem
  \begin{align*}
    \max_{b\leq r\leq b+1} |v(r)|^2 r 
    \leq \Phi(b) \int_b^{b+1} |v'(r)|^2 r \, dr,
  \end{align*}
  satisfies $\max_{b\leq r\leq b+1} |v_*(r)|^2 r = |v_*(b+1)|^2 (b+1)$.
  The proof of Proposition~\ref{propositionA} would be
  simplified if we could prove this a priori.
\end{remark}


\section{A class of half-line operators}

Let $\sigma\in\R$. We consider the quadratic form 
\begin{align}
h_\sigma [u] 
:= \int_{\R_+} \left( r^\sigma |u'(r)|^2 - \frac{(\sigma -1)^2
  |u(r)|^2}{4 r^{2 - \sigma}} \right) \, dr 
\end{align}
defined on $C_0^\infty(\R_+)$ if $\sigma\leq 1$ and on
$C_0^\infty(\overline{\R_+})$ if $\sigma> 1$. The generalized Hardy
inequality \eqref{eq:hardygeneral} implies that the forms $h_\sigma$
are non-negative on their respective domains. (We shall essentially
reprove this in the proof of Theorem \ref{lthsigma}.) Moreover,
they are closable in $L_2(\R_+)$ and thus generate self-adjoint
operators $H_\sigma$. Note that the operator $H_0$ coincides with the
operator $- \frac {d^2}{dr^2} - \frac 1{4r^2}$ treated in previous
sections. 

Our main result in this section are Lieb-Thirring inequalities on the
moments of negative eigenvalues of the Schr\"odinger-type operator
$H_\sigma-V$. 

\begin{theorem}\label{lthsigma}
Let $\gamma > 0$. 
Assume that either
\begin{equation*}
\sigma > 2,\ \gamma>0,\ \alpha \leq - \frac \sigma 2
\quad\textit{satisfy}\quad \gamma - \frac{1+\alpha}{\sigma-2} \geq 1
\end{equation*}
or otherwise that
\begin{equation*}
\sigma < 2,\ \gamma>0,\ \alpha \geq - \frac \sigma 2
\quad\textit{satisfy}\quad \gamma + \frac{1+\alpha}{2-\sigma} \geq 1.
\end{equation*}
Then
\begin{align} \label{eq:lthsigma}
\tr (H_\sigma - V)_-^\gamma 
\leq \left|\frac{2}{2 - \sigma}\right|^{\frac{2\alpha+\sigma}{2-\sigma}}
C_{\gamma,\frac{2\alpha+\sigma}{2-\sigma}} 
\int_{\R_+} V(r)_+^{\gamma + \frac {1+\alpha}{2 - \sigma}} r^\alpha \, dr,
\end{align}
where the constant $C_{\gamma, \alpha}$ is given in Theorem \ref{maintheorem}.

Finally, let $\sigma=2$ and assume that either $\gamma \geq 1/2$ and
$\alpha=0$ or otherwise that $\gamma>0$ and $\alpha>0$ with $\gamma +
\frac{1+\alpha}{2}\geq 1+\alpha$. Then 
\begin{align} \label{eq:lthsigma2}
\tr (H_2 - V)_-^\gamma \leq C_{\gamma,\alpha}^{EK} \int_{\R_+}
V(r)_+^{\gamma + \frac{1+\alpha}2} |\log r|^\alpha r^{-1} \, dr,
\end{align}
where the constant $C_{\gamma,0}^{EK}=L_\gamma$ is given in \eqref{eq:lt} and $C_{\gamma,\alpha}^{EK}$ in \eqref{eq:ek}.
\end{theorem}

\begin{proof}
We shall first consider the case $\sigma\neq 2$ and introduce the unitary operators $U_\sigma$ in $L_2(\R_+)$,
\begin{align*}
	(U_\sigma u)(r) := \left| \frac {2-\sigma}{2} \right|^{1/2}
  r^{-\frac{\sigma}{4}} u(r^\frac{2-\sigma}{2}).
\end{align*}
We note that the adjoint operators $U_\sigma^*$ map $C_0^\infty(\R_+)$
into itself and, for $\sigma>1$, map  
$C_0^\infty(\overline{\R_+})$ into the form domain of $H_0$. Moreover, one easily verifies that
\begin{equation*}
	h_\sigma[u] = \left(\frac{2-\sigma}2\right)^2 h_0[U_\sigma^* u].
\end{equation*}
This relation, proved initially for $u\in C_0^\infty(\R_+)$ if $\sigma\leq 1$ -- or for $u\in C_0^\infty(\overline{\R_+})$ if $\sigma> 1$ -- extends to the closure of $h_\sigma$ and implies that
\begin{align*}
	U_\sigma^* H_\sigma U_\sigma = \left( \frac {2 - \sigma}{2} \right)^2 H_0.
\end{align*}
For given $V$ we define $V_\sigma(r):=  \left(\frac2{2-\sigma}\right)^2 V\left( r^\frac 2{2-\sigma}\right)$ and find that
\begin{align*}
	\tr (H_\sigma - V)_-^\gamma & 
	= \left| \frac {2 - \sigma}{2} \right|^{2 \gamma} \tr \left( H_0 - V_\sigma \right)_-^\gamma.
\end{align*}
Hence \eqref{eq:lthsigma} follows from our main result, Theorem \ref{maintheorem}.

In the case $\sigma = 2$ we define the unitary operator $U_2 : L_2(\R) \to L_2(\R_+)$ by $(U_2 u)(r) := r^{-\frac 12} u(\log r)$. Similarly as above one checks that
\begin{align*}
U_2^* H_2 U_2 = - \frac {d^2}{dx^2}
\end{align*}
and hence with $V_2(x):= V(e^x)$,
\begin{align*}
\tr_{L_2(\R_+)} (H_2 - V)_-^\gamma 
= \tr_{L_2(\R)} \left( -\frac{d^2}{dx^2} - V_2 \right)_-^\gamma.
\end{align*}
Inequality \eqref{eq:lthsigma2} now follows from \eqref{eq:lt} and \eqref{eq:ek}.
\end{proof}

\begin{remark}
Note that the sharp constants $C_{\gamma,0}^{EK}=L_\gamma$ in \eqref{eq:lthsigma2} with $\alpha=0$ are known if $\gamma
=1/2$ or $\gamma\geq 3/2$.
\end{remark}

\begin{remark}
The method of the previous proof allows one to obtain rather complete information on the behavior of weakly coupled eigenvalues. Assume for simplicity that $V$ is bounded with compact support in $\R_+$.\footnote{It will be sufficient that $\int_{\R_+} |V(r)|^{1+\delta} r^{1-\sigma} \, dr + \int_{\R_+} |V(r)| (1 + r^{(2-\sigma)\delta/2}) r^{1-\sigma} \, dr < \infty$ for some $\delta > 0$ if $\sigma=2$ and that $\int_{\R_+} |V(r)| (1+r)\frac {dr}r <\infty$ if $\sigma=2$.}
Then $H_\sigma - \beta V$ has a negative eigenvalue for all $\beta > 0$ if and only if $V \not\equiv 0$ and $\int_{\R_+} V(r) r^{1 - \sigma}\, dr \geq 0$. In the case $\int_{\R_+} V(r) r^{1-\sigma} \,dr > 0$ there is a unique eigenvalue $\lambda(\beta)$ for all sufficiently small $\beta$, and one has
\begin{align*}
& \lim_{\beta \to 0} \beta^{-1} |\log |\lambda(\beta)||^{-1} 
= \frac1 {2 - \sigma} \int_{\R_+} V(r) r^{1 - \sigma} \, dr,
\quad \sigma\neq 2, \\
& \lim_{\beta \to 0} \beta^{-2} \lambda(\beta) 
= - \frac 14 \left( \int_{\R_+} V(r) \, \frac {dr}{r} \right)^2,
\quad \sigma=2.
\end{align*}
Like in the proof of Proposition \ref{lthsigma} this assertion is reduced to the operators $H_0$ and $-\frac{d^2}{dx^2}$, for which the assertion is known, see \cite{Si1} and \cite{EkFr}. One can treat the case where $\int_{\R_+} V(r) r^{1-\sigma} \,dr = 0$ in a similar manner, but we omit the details.
\end{remark}



\bibliographystyle{amsalpha}

\end{document}